\documentclass[psamsfont]{amsart}

\usepackage{amsmath,amssymb}

\usepackage{hyperref}

\newtheorem{theo}{Theorem}[section]
\newtheorem{claim}[theo]{Claim}
\newtheorem{claim2}[theo]{Claim}
\newtheorem{defi}[theo]{Definition}
\newtheorem{lemm}[theo]{Lemma}
\newtheorem{coro}[theo]{Corollary}
\newtheorem*{exam}{Example}

\newtheorem{rem}[theo]{Remark}

\newenvironment{rema}{\begin{rem}\rm}{\end{rem}}

\usepackage{xcolor,marvosym}
\let\oldmarginpar\marginpar
\renewcommand\marginpar[1]{\-\oldmarginpar[\raggedleft #1]%
{\raggedright #1}}

\newcommand{\defin}[1]{\begin{color}{black}{\textsf{#1}}\end{color}}

\usepackage{graphicx}

\usepackage[all]{xy}

\def\del{\partial} 
\def\bb#1{\mathbb{#1}} \def\m#1{\mathcal{#1}}
\def\id{\operatorname{id}}
\def\im{\operatorname{im}}
\def\ra{\rightarrow} \def\lra{\longrightarrow} \def\co{\colon\thinspace}

\def\hop{{\noindent}}

\def\ham#1{\mathrm{Ham}(#1)}
\def\PSS{{\rm PSS}}
\def\pop{\star_{\footnotesize{\rm pp}}}
\def\seidel{{\bf S}}
\def\nat{{\bf n}}
\def\comp{{\bf c}}

\begin{document}

\title[The Seidel morphism of products]{The Seidel morphism of Cartesian products}
\author{R\'emi Leclercq}
\address{R.L.: Max-Planck-Institut f\"ur Mathematik in den Naturwissenschaften, Inselstra\ss e 22, 04103 Leipzig, Deutschland}
\email{leclercq@mis.mpg.de}

\begin{abstract}
We prove that the Seidel morphism of $(M\times M',\omega\oplus\omega')$ is naturally related to the Seidel morphisms of $(M,\omega)$ and $(M',\omega')$, when these manifolds are monotone. We deduce that any homotopy class of loops of Hamiltonian diffeomorphisms of one component, with non-trivial image via Seidel's morphism, leads to an injection of the fundamental group of the group of Hamiltonian diffeomorphisms of the other component into the fundamental group of the group of Hamiltonian diffeomorphisms of the product. This result was inspired by and extends results obtained by Pedroza \cite{Pedroza08}.
\end{abstract}

\maketitle

\section{Introduction}

All the symplectic manifolds we consider are closed. Let $(M,\omega)$ be a $2n$--dimensional closed symplectic manifold, and let $c_1$ denote the first Chern class of its tangent bundle, $c_1(TM,\omega)$. The \defin{minimal Chern number} of $(M,\omega)$ is then defined as the infimum
\begin{align*}
N = \inf \{ k>0 \,|\; \exists A\in\pi_2(M),\, c_1(A)=k \}
\end{align*} 
and we put $N=\infty$ if the latter set is empty. Notice that, when $N$ is finite, $c_1(\pi_2(M))=N\bb Z$. The symplectic manifold $(M,\omega)$ is \defin{strongly semi-positive} if (at least) one of the following conditions holds:
\begin{itemize}
\item[(a)] there exists $\lambda\geq 0$, such that for all $A\in\pi_2(M)$, $\omega(A)=\lambda c_1(A)$,
\item[(b)] $c_1$ vanishes on $\pi_2(M)$,
\item[(c)] the minimal Chern number satisfies $N\geq n-1$.
\end{itemize}
Under this assumption, Seidel introduced \cite{Seidel97} a group morphism:
\begin{align*}
q_M\co \tilde\pi_1(\ham{M,\omega})\lra QH_*(M,\omega)^\times,
\end{align*}
where $QH_*(M,\omega)^\times$ denotes the group of invertible elements of $QH_*(M,\omega)$, the quantum homology of $(M,\omega)$. We recall that the identity of the group $QH_*(M,\omega)^\times$ is the fundamental class of $M$, which is denoted $[M]$.

As usual, $\ham{M,\omega}$ denotes the group of Hamiltonian diffeomorphisms of $(M,\omega)$ and $\tilde{\pi}_1(\ham{M,\omega})$ is a covering of $\pi_1(\ham{M,\omega})$ which will be defined below. The inclusions of $\ham{M,\omega}$ and $\ham{M',\omega'}$ in $\ham{M\times M',\omega\oplus \omega'}$ induce a map between the respective fundamental groups: $([g],[g'])\mapsto [g,g']$, where $[g,g']$ stands for $[(g,g')]$, the homotopy class of the loop $(g,g')$. The extension of this map to the coverings $\tilde\pi_1$ is straightforward. We denote it by
\begin{align*}
i\co \tilde{\pi}_1(\ham{M,\omega}) \times \tilde{\pi}_1(\ham{M',\omega'})\lra \tilde{\pi}_1(\ham{M\times M',\omega\oplus \omega'}).
\end{align*}

We also denote by
\begin{align*}
\kappa_Q\co QH_*(M,\omega)\otimes QH_*(M',\omega') \lra QH_*(M\times M',\omega\oplus\omega')
\end{align*}
the inclusion given by K\"unneth formula and the compatibility of the Novikov rings with the Cartesian product (see \S \ref{sec:basics} for definitions).\\

Let $(M,\omega)$ and $(M',\omega')$ be strongly semi-positive symplectic manifolds and let $[\bar{g}]\in\tilde{\pi}_1(\ham{M,\omega})$ and $[\bar{g}']\in\tilde{\pi}_1(\ham{M',\omega'})$ (which respectively lift the homotopy classes of loops $g\subset\ham{M,\omega}$ and $g'\subset\ham{M',\omega'}$). When $(M\times M',\omega\oplus\omega')$ is strongly semi-positive (this is not necessarily the case, see discussion in Remark \ref{rema:MxM satisfies W+}), one can, on one hand, compute the images of $[\bar{g}]$ and $[\bar{g}']$ via the respective Seidel's morphisms and then see the result as an element in $QH_*(M\times M',\omega\oplus\omega')^\times$ via $\kappa_Q$. On the other hand, one can compute the image of $i([\bar{g}],[\bar{g}'])$, via the Seidel morphism of the product.

In this note, we prove that both computations coincide when the involved manifolds are monotone. A symplectic manifold is called \defin{monotone} if it satisfies condition (a) above, with $\lambda> 0$. Notice that, if $(M,\omega)$ and $(M',\omega')$ are monotone with constants $\lambda$ and $\lambda'$, the product $(M\times M',\omega\oplus\frac{\lambda}{\lambda'}\omega')$ is monotone. In particular, if $\lambda=\lambda'$, $(M\times M',\omega\oplus\omega')$ is monotone.

\begin{theo}\label{theo:main theo}
Let $(M,\omega)$ and $(M',\omega')$ be closed monotone symplectic manifolds (with identical monotonicity constants), and let $[\bar{g}]\in \tilde{\pi}_1(\ham{M,\omega})$ and $[\bar{g}'] \in \tilde{\pi}_1(\ham{M',\omega'})$, then
\begin{align*}
  q_{M\times M'}\big(i([\bar{g}],[\bar{g}'])\big) = \kappa_Q\big(q_M([\bar{g}])\otimes q_{M'}([\bar{g}'])\big).
\end{align*}
\end{theo}

Notice that this statement is conditioned by coherent choices of the three involved Novikov rings (see definition in \S \ref{sec:basics} below), since several slightly different versions are commonly used.

\begin{rema}\label{rema:on the monotonicity restriction}
The monotonicity assumption ensures that there exists no non-constant pseudo-holomorphic sphere of first Chern number $0$. This property is only used in the proof of Lemma \ref{lemm:regularity of J+J'}, which states that a particular choice of almost complex structures is regular enough to compute Seidel's morphism. Thus, all the results of this note hold under the weaker assumption that both manifolds \textit{and} their product are ``strongly semi-positive, without non-constant pseudo-holomorphic spheres with first Chern number $0$''.

For strongly semi-positive manifolds admitting such spheres, the theorem is more difficult to prove but most probably holds, see Remark \ref{rema:SSP non monotone} (we do mean that it holds even without the use of virtual techniques, see Remark \ref{rema:virtual generalization}).
\end{rema}

Even though the morphism $q$ is interesting in itself, one usually looks for information concerning $\pi_1(\ham{M,\omega})$ (rather than $\tilde{\pi}_1(\ham{M,\omega})$). Now, the definition of quantum homology relies on a Novikov ring built from $\Gamma$, a quotient of $\pi_2(M)$. Moreover, $\Gamma$ can be seen as a subgroup of the group of invertible elements of quantum homology via the map $\tau$, defined by $\tau(\gamma)=[M]\otimes \gamma$ for all $\gamma\in\Gamma$. Seidel's morphism then induces a morphism $\bar{q}$ defined by the commutativity of the diagram:
  \begin{align}\label{eq:diagram defining qbar}
    \begin{split}
      \xymatrix{\relax
        \tilde{\pi}_1(\ham{M,\omega}) \ar[d]_{q}\ar[r] & \pi_1(\ham{M,\omega})\ar[d]^{\bar{q}}\\
        QH_{*}(M)^\times \ar[r] & QH_{*}(M)^\times /\tau(\Gamma)
      }
    \end{split}
  \end{align}

The main consequence of Theorem \ref{theo:main theo} can be stated in terms of $\bar{q}$.
\begin{coro}\label{coro:main corollary}
  Let $(M,\omega)$, $(M',\omega')$ be as in the theorem. Let $g$ and $g'$ be loops of Hamiltonian diffeomorphisms of respectively $(M,\omega)$ and $(M',\omega')$. If $[g,g']$ is trivial in $\pi_1(\ham{M\times M',\omega\oplus \omega'})$, both loops $g$ and $g'$ are mapped to the identity via Seidel's morphism, namely, $\bar{q}_M([g])=[M]$ and $\bar{q}_{M'}([g'])=[M']$.
\end{coro}

We will deduce this corollary from the theorem, at the end of \S \ref{subsec:def Seidel morphism}, as soon as all the involved objects have been defined.

\begin{rema}\label{rema:virtual generalization}
Notice that, in order to get (the same) results \textit{directly on the fundamental groups}, one could also adopt the approach of Lalonde, McDuff and Polterovich \cite{LalondeMcDuffPolterovich99}. Furthermore, following McDuff \cite{McDuff00},  one could address these questions by means of virtual techniques. This would probably provide a proof of the results contained in this note in great generality. We thank Shengda Hu for pointing this fact out to us.
\end{rema}

From Corollary \ref{coro:main corollary}, it is easy to derive the following properties.

\begin{coro}\label{coro:M and M'}
  Let $(M,\omega)$ be monotone and let $[g]\in\pi_1(\ham{M,\omega})$ such that $\bar{q}_M(g)\neq [M]$. Then for any monotone symplectic manifold $(M',\omega')$ (with the same monotonicity constant), the map
  \begin{align*}
  \iota{[g]}\co\pi_1(\ham{M',\omega'}) \lra \pi_1(\ham{M\times M',\omega\oplus\omega'}),
  \end{align*}
  defined by $\iota{[g]}([g'])=[g,g']$, is injective. Moreover, if $\bar{q}_M([g_1])\neq \bar{q}_M([g_2])$, then
  \begin{align*}
  \iota{[g_1]}(\pi_1(\ham{M',\omega'}))\cap \iota{[g_2]}(\pi_1(\ham{M',\omega'}))=\emptyset.
  \end{align*}
\end{coro}

Let us emphasize the fact that these maps $\iota{[g]}$ are maps between sets (and not group morphisms).

\begin{coro}\label{coro:M and M}
  Let $(M,\omega)$ be a monotone symplectic manifold. If Seidel's morphism $\bar{q}_M$ is injective, then the map
  \begin{align*}
        \pi_1(\ham{M,\omega})\times \pi_1(\ham{M,\omega}) \lra \pi_1(\ham{M\times M,\omega\oplus\omega})
  \end{align*}
  induced by the inclusion is injective.
\end{coro}

\begin{exam}
The complex projective spaces $\bb C\bb P^1$ and $\bb C\bb P^2$, endowed with the Fubini--Study symplectic form, are symplectic manifolds with injective Seidel's morphism and, more generally, the whole family of complex projective spaces provides interesting examples.

More precisely, let $\omega_\mathrm{st}$ be the Fubini--Study symplectic form on $\bb C\bb P^m$, normalized in such a way that $(\bb C\bb P^{m},\omega_\mathrm{st})$ is monotone, with monotonicity constant $1/(m+1)$. Seidel \cite{Seidel97} proved that there exists an element of order $m+1$ in $\pi_1(\ham{\bb C\bb P^m,\omega_\mathrm{st}})$\footnote{Actually, Seidel proved the existence of such an element for any complex Grassmannian $\mathrm{Gr}_k(\bb C^{m+1})$, $k\geq 1$, endowed with the symplectic form induced by Pl\"ucker's embedding.}. This explicit element comes from the action of $\mathrm{U}(m+1)$ on $\bb C^{m+1}$; we denote it by $\alpha_m$.

In order to obtain monotone products, we consider a multiple of the standard symplectic form, namely, we endow $\bb C\bb P^{m}$ with $\omega_m=(m+1)\,\omega_\mathrm{st}$. We also denote the element of order $m+1$ of $\pi_1(\ham{\bb C\bb P^m,\omega_m})$ by $\alpha_m$.

From Corollary \ref{coro:M and M'} we deduce the following properties.
\begin{itemize}
\item[(1)] For $1\leq l\leq n$, the inclusions $\iota{[\alpha_n^l]}$, given by $\iota{[\alpha_n^l]}([g'])= [\alpha_n^l,g']$, lead to $n$ distinct copies of $\pi_1(\ham{\bb C\bb P^{n'},\omega_{n'}})$ in $\pi_1(\ham{\bb C\bb P^{n}\times \bb C\bb P^{n'},\omega_{n}\oplus\omega_{n'}})$.
\item[(1')] Similarly, $\pi_1(\ham{\bb C\bb P^{n}\times \bb C\bb P^{n'},\omega_{n}\oplus\omega_{n'}})$ also contains $n'$ distinct copies of $\pi_1(\ham{\bb C\bb P^{n},\omega_{n}})$, given by $\iota'{[\alpha_{n'}^{l'}]}([g])= [g,\alpha_{n'}^{l'}]$, for $1\leq l'\leq n'$.
\item[(2)] These injections intersect pairwise in a unique point, namely, $\iota{[\alpha_n^l]}([g'])=\iota'{[\alpha_{n'}^{l'}]}([g])$ if and only if $g=\alpha_n^l$ and $g'=\alpha_{n'}^{l'}$.
\item[(3)] Finally, the elements $(\alpha_n,\id_{\bb C\bb P^{n'}})$ and $(\id_{\bb C\bb P^{n}},\alpha_{n'})$ are of respective orders $n+1$ and $n'+1$ (and $\pi_1(\ham{\bb C\bb P^{n}\times \bb C\bb P^{n'},\omega_{n}\oplus\omega_{n'}})$ contains subgroups isomorphic to $\bb Z_{n+1}$ and $\bb Z_{n'+1}$). 
\end{itemize}
Notice that when $n=n'$ the same statements hold with $\omega_\mathrm{st}$. 

Finally, it is well-known that $\pi_1(\ham{\bb C\bb P^1,\omega_\mathrm{st}})=\bb Z_2$ and $\pi_1(\ham{\bb C\bb P^2,\omega_\mathrm{st}})=\bb Z_3$. (In the latter case, Gromov \cite{Gromov85} proved that the group of symplectomorphisms of $(\bb C\bb P^2,\omega_\mathrm{st})$ contracts onto the group of isometries, which is isomorphic to $PU(3)$. Now, $H^1_\mathrm{dR}(\bb C\bb P^2)=0$ and thus, the symplectomorphisms of $\bb C\bb P^2$ -- which are isotopic to the identity -- are Hamiltonian diffeomorphisms.) Since Seidel's morphism of $(\bb C\bb P^n,\omega_\mathrm{st})$ detects an element of order $n+1$, it follows that it is injective in these particular cases. Corollary \ref{coro:M and M} allows us to conclude that the obvious mapping
\begin{align*}
\pi_1(\ham{\bb C\bb P^{n},\omega_\mathrm{st}})\times \pi_1(\ham{\bb C\bb P^{n},\omega_\mathrm{st}}) \longrightarrow \pi_1(\ham{\bb C\bb P^{n}\times \bb C\bb P^{n},\omega_\mathrm{st}\oplus\omega_\mathrm{st}})
\end{align*}
is injective, when $n=1$ and $n=2$.
\end{exam}

\begin{rema}
Another application of Theorem \ref{theo:main theo} comes from recent work of Hu and Lalonde \cite{HuLalonde05}. Indeed, they introduced a relative version (that is, defined with respect to a Lagrangian $L$) of Seidel's morphism and they proved that it is related to the Seidel morphism of the ambient manifold $(W,\Omega)$ via a map defined by Albers \cite{Albers05} (under suitable assumptions).

Let $(W,\Omega)=(M\times M,\omega\oplus(-\omega))$ and $L$ be the graph of a Hamiltonian diffeomorphism of $(M,\omega)$. Combining Theorem \ref{theo:main theo} with the morphism introduced by Biran, Polterovich and Salamon \cite{BiranPolterovichSalamon03} allows us to compare the relative Seidel morphism associated to $L$, not only to the absolute Seidel morphism associated to $(W,\Omega)$ but also to the one associated to $(M,\omega)$ (at least for ``split'' loops).
\end{rema}

\begin{rema}\label{rema:MxM satisfies W+}
As mentioned above, being strongly semi-positive is not \emph{a priori} compatible with the Cartesian product. Nevertheless, as for monotone symplectic manifolds, the product of certain manifolds is automatically strongly semi-positive. Let, for instance, $(M,\omega)$ and $(M',\omega')$ both satisfy the (sub-)condition
\begin{itemize}
\item[(a)] with constants $\lambda$ and $\lambda'$. If $\lambda=\lambda'=0$, the product also satisfies this condition. (In that case, there is no non-constant pseudo-holomorphic sphere since such a sphere has positive symplectic area. Thus, Seidel's morphism is trivial for these manifolds and Theorem \ref{theo:main theo} is trivially satisfied.)
\item[(b)] $c_1$ vanishes on $\pi_2(M)$, $c'_1$ vanishes on $\pi_2(M')$. Thus the first Chern class of the tangent bundle of $M\times M'$ vanishes on $\pi_2(M\times M')$ and the product is strongly semi-positive.
\item[(c)] the minimal Chern numbers $N$ and $N'$ satisfy $N\geq n-1$ and $N'\geq n'-1$. The minimal Chern number of the product being the greatest common divisor of $N$ and $N'$, for the product to satisfy sub-condition (c), the gcd of $N$ and $N'$ has to satisfy $\mathrm{gcd}(N,N')\geq n+n'-1$. Notice for example that if $N\geq 2n-1$, then $(M\times M,\omega\oplus\omega)$ is strongly semi-positive.
\end{itemize}
\end{rema}

There is also another remarkable particular case: when $\pi_2(M')=0$, then of course if $(M,\omega)$ is monotone (respectively, satisfies (a), (b) or (c) with $N\geq n+n'-1$), the product is monotone (respectively, strongly semi-positive). Thus Theorem \ref{theo:main theo} extends results obtained by Pedroza. Indeed, \cite[Theorem 1.1]{Pedroza08} is given by our theorem, in the case where $\pi_2(M')=0$, $\bar{g}'=\id$, and \cite[Theorem 1.3]{Pedroza08} corresponds to the case where $M=M'$, $\pi_2(M)=0$ and $\bar{g}'=\bar{g}$.

Hence, our extension is two-fold. Firstly, we do not have to restrict ourselves to special loops among the split loops of Hamiltonian diffeomorphisms of the product, the homotopy class of $(\bar{g},\bar{g}')$ being detected by Seidel's morphism as soon as $[\bar{g}]$ or $[\bar{g}']$ is detected. Secondly, concerning the latter result, we emphasize the fact that we \emph{do not require any type of asphericity condition} such as $\pi_2(M)$ trivial. This is important since, for $(M',\omega')=(M,\omega)$ with $\pi_2(M)=0$, there is no non-trivial (pseudo-holomorphic) sphere and the involved Seidel morphisms are trivial. Thus, the statement of Theorem \ref{theo:main theo} is, in that very particular case, trivially satisfied.

Another noteworthy difference between this note and \cite{Pedroza08} is the approach to Seidel's morphism which we consider. Pedroza approaches the question via the point of view of Hamiltonian fibrations, we use the representation approach (in terms of automorphisms of Floer homology).\\

In the next section we recall the definitions of quantum homology (\S \ref{subsec:def QH}) and of Seidel's morphism (\S \ref{subsec:def Seidel morphism}). This allows us to prove Corollary \ref{coro:main corollary} from Theorem \ref{theo:main theo}. Then we recall the construction of Floer homology and the representation viewpoint on Seidel's morphism (\S \ref{subsec:def seidel representation}). Afterwards, we prove Theorem \ref{theo:main theo} (\S \ref{sec:proof of theo}) up to a claim concerning the regularity of a particular choice of parameters (required to compute Seidel's morphism). Finally, we justify the claim (\S \ref{section:regularity of split pairs}).


\section{Making things precise}\label{sec:basics}

\subsection{The group $\Gamma$ and the morphism $\kappa_Q$}\label{subsec:def QH}

Following Seidel, we define $\Gamma_M$ as the group of equivalence classes of elements in $\pi_2(M)$ under the equivalence relation $A\sim B$ if $\omega(A)=\omega(B)$ and $c_1(A)=c_1(B)$. Notice that the obvious map
\begin{align*}
\pi_\Gamma\co\Gamma_{M}\times \Gamma_{M'} &\lra \Gamma_{M\times M'}\\
([A],[A'])&\longmapsto [A,A']
\end{align*}
is well-defined and surjective but in general not injective. Its kernel consists of pairs $([A],[A'])$ such that $\omega(A)+\omega'(A')=0$ and $c_1(A)+c'_1(A')=0$.

We recall that the (small) quantum homology is the $\Lambda_M$--module given as the tensor product $H_*(M,\bb Z_2)\otimes_{\bb Z_2} \Lambda_M$, where $\Lambda_M$ is the Novikov ring defined as the group of formal sums $\sum m_\gamma \cdot\gamma$, with $\gamma\in\Gamma_M$, $m_\gamma\in \bb Z_2$ and satisfying the finiteness condition:
\begin{align*}
\forall C\in \bb R,\; \#\{m_\gamma |\, m_\gamma\neq 0,\, \omega(\gamma)\leq C\} < \infty.
\end{align*}
Since an element of the type $[M]\otimes \gamma$ is invertible in $QH_*(M,\omega)$, the formula $\tau(\gamma)= [M]\otimes \gamma$ defines a morphism $\tau\co \Gamma_M \ra QH_*(M,\omega)^\times$. Since this morphism is injective, $\Gamma_M\simeq\tau(\Gamma_M)$ can be seen as a subgroup of $QH_*(M,\omega)^\times$.\\

Now let $\kappa\co H_*(M,\bb Z_2)\otimes_{\bb Z_2} H_*(M',\bb Z_2)\lra H_*(M\times M',\bb Z_2)$ be the inclusion given by K\"unneth formula (due to the fact that we use the field $\bb Z_2$ for coefficients, this is actually an isomorphism). We define $\kappa_Q$ by the formula
\begin{align*}
\kappa_Q\co QH_*(M,\omega)\otimes_{\bb Z_2} QH_*(M',\omega') &\lra QH_*(M\times M',\omega\oplus\omega')\\
\big((\alpha\otimes\gamma)\otimes(\alpha'\otimes\gamma')\big) &\longmapsto \kappa(\alpha\otimes\alpha')\otimes \pi_\Gamma(\gamma,\gamma')
\end{align*}
for simple tensors and extend it by linearity for general quantum elements. By definition of $\kappa_Q$ and injectivity of $\kappa$, we deduce the following lemma.
\begin{lemm}\label{lemm:injectivity kappa's}
Let $\kappa_Q (q\otimes q')=[M\times M']\otimes\gamma_\times$ for some $\gamma_\times\in \Gamma_{M\times M'}$. There exist $\lambda\in\Lambda_M$ and $\lambda'\in\Lambda_{M'}$ such that $q=[M]\otimes\lambda$ and $q'=[M']\otimes\lambda'$.
\end{lemm}

\subsection{Seidel's morphism and the proof of Corollary \ref{coro:main corollary}}\label{subsec:def Seidel morphism}

Let $(M,\omega)$ be a monotone symplectic manifold. Following Seidel's notation, $G$ denotes the set of smooth loops of Hamiltonian diffeomorphisms (based at the identity). Then $\pi_0(G)\simeq \pi_1(\ham{M,\omega})$. Let $\m LM$ be the set of free, smooth, contractible loops of $M$. $\widetilde{\m LM}$ is the set of equivalence classes of pairs $(v,x)\in C^\infty(D^2,M)\times \m LM$ such that $v$ coincides with $x$ on the boundary $\del D^2$, under the equivalence relation
\begin{align*}
(v,x)\sim(v',x') \Leftrightarrow x=x' \;\mathrm{and}\; \omega(v\#\overline{v'})=0, c_1(v\#\overline{v'})=0.
\end{align*}
Here $\overline{v'}$ is $v'$ considered with the opposite orientation and $v\#\overline{v'}$ the sphere obtained by gluing the two discs along their common boundary. There is an action of $G$ on $\m LM$, given by $g\cdot x = [t\mapsto g_t(x(t))]$, which lifts to $\widetilde{\m LM}$. We define $\widetilde{G}$ as the subset of $G\times \mathrm{Homeo}(\widetilde{\m LM})$ consisting of pairs $(g,\tilde{g})$ such that $\tilde{g}$ is a lift of $g$ (that is, $\tilde{g}(v,x)=(v',g\cdot x)$ for any $(v,x)\in \widetilde{\m LM}$).

$\pi_0(\widetilde{G})$ is the covering of $\pi_1(\ham{M,\omega})$ which was denoted $\tilde\pi_1(\ham{M,\omega})$ above. Following Witten \cite{Witten88}, Seidel introduced in \cite{Seidel97} the morphism $q$ by considering Hamiltonian fibre bundles over $S^2$, with fibre $(M,\omega)$. Roughly, since a fibre bundle over a disc is trivial, it is easy to see that such a fibre bundle over $S^2$ corresponds to the choice of a loop $g$ of Hamiltonian diffeomorphisms (based at the identity). Seidel then derived invariants from the pseudo-holomorphic sections of these bundles by comparing them to some chosen equivalence class of sections (given by the choice of a lift of $g$).

\begin{rema}
Now that $\tilde\pi_1(\ham{M,\omega})$ has been made precise, we will also use the following obvious notation $i([g,\tilde{g}],[g',\tilde{g}']) = [g,g';\tilde{g},\tilde{g}']$ to denote the morphism $i$ defined above.
\end{rema}

We can now deduce Corollary \ref{coro:main corollary} from Theorem \ref{theo:main theo}.

\begin{proof}[Proof of Corollary \ref{coro:main corollary}]
Let $g$ and $g'$ be loops of Hamiltonian diffeomorphisms such that $[g,g']$ is trivial in $\pi_1(\ham{M\times M',\omega\oplus \omega'})$. Then $\bar{q}([g,g'])$ is the identity of $QH_*(M\times M',\omega\oplus\omega')^\times / \tau(\Gamma_{M\times M'})$. By definition of $\bar{q}$ (that is, by commutativity of Diagram (\ref{eq:diagram defining qbar})), this amounts to the fact that for any lift $(g,g';\tilde{g},\tilde{g}')$ of $(g,g')$, $q([g,g';\tilde{g},\tilde{g}'])\in\tau(\Gamma_{M\times M'})$. We fix $\tilde{g}$ and $\tilde{g}'$ respective lifts of $g$ and $g'$.

We know that there exists $\gamma_\times\in \Gamma_{M\times M'}$ such that $q([g,g';\tilde{g},\tilde{g}'])=[M\times M']\otimes \gamma_\times$. By Theorem \ref{theo:main theo}, this gives that $\kappa_Q\big(q_M([g,\tilde{g}])\otimes q_{M'}([g',\tilde{g}'])\big) =[M\times M']\otimes \gamma_\times$.

By Lemma \ref{lemm:injectivity kappa's}, this implies the existence of $\lambda\in\Lambda_M$ and $\lambda'\in\Lambda_{M'}$ such that $q_M([g,\tilde{g}])=[M]\otimes\lambda$ and $q_{M'}([g',\tilde{g}'])=[M']\otimes\lambda'$. Since
\begin{align*}
[M]\otimes\lambda=[M]\otimes\sum m_{\gamma}\gamma=\sum m_\gamma ([M]\otimes \gamma),
\end{align*}
$q_M([g,\tilde{g}])$ is mapped to the identity element in $QH_{*}(M)^\times /\tau(\Gamma_M)$ (and similarly for $q_{M'}([g',\tilde{g}'])$); by definition of $\bar{q}$, so are $\bar{q}_M([g])$ and $\bar{q}_{M'}([g'])$.
\end{proof}

\subsection{An alternate description: Seidel's representation}\label{subsec:def seidel representation}

As noticed by Seidel, there is an alternate description of $q$ as a representation of $\tilde{\pi}_1(\ham{M,\omega})$ in terms of automorphisms of Floer homology. \\

We briefly recall the definition of the Floer homology of a closed (monotone) symplectic manifold $(M,\omega)$. Let $H$ be a Hamiltonian function on $M$. The action functional is defined on $\widetilde{\m LM}$ by the formula
\begin{align*}
\m A_H([v,x]) = -\int v^*\omega +\int H_t(x(t))dt.
\end{align*}
The set of critical points of $\m A_H$, $\mathrm{Crit}(\m A_H)$, consists of equivalence classes $[v,x]$ where $x$ is a contractible periodic orbit of $X_H$, the Hamiltonian vector field generated by $H$. We recall that $\mathrm{Crit}(\m A_H)$ is graded via the Conley--Zehnder index (see \cite{ConleyZehnder83} for definition and, for example, \cite{Salamon99} for another comprehensive description). 

We now pick a regular almost complex structure $J$ on $TM$ such that the pair $(H,J)$ is regular (which means, in particular, that the critical points of $\m A_H$ are non-degenerate, see \S \ref{section:regularity of split pairs} for precise definitions). A generic choice of pair is regular and for such a choice, one can define the Floer complex $(CF_*(H),\del_{(H,J)})$ where $CF_k(H)$ is the group of formal sums $\sum_c m_c\cdot c$ where the $c$'s are critical points of $\m A_H$ of index $k$, $c\in \mathrm{Crit}_k(\m A_H)$. The coefficients $m_c$ are elements of $\bb Z_2$ and satisfy the finiteness condition
\begin{align}\label{eq:finiteness condition FH}
\#\{ m_c|\, m_c\neq 0 \;\mathrm{and}\; \m A_H(c)\geq C \}<\infty
\end{align}
for all real numbers $C$.

The differential $\del_{(H,J)}\co \mathrm{Crit}_k(\m A_H)\ra \mathrm{Crit}_{k-1}(\m A_H)$ is defined by the formula
\begin{align*}
\del_{(H,J)}(c) = \sum_{c'\in \mathrm{Crit}_{k-1}(\m A_H)} \#_2[\m M(c,c';H,J)/\bb R]\cdot c'
\end{align*}
where $\#_2\m M(c,c';H,J)$ is the cardinality (mod $2$) of the set of curves $u\co S^1\times \bb R\ra M$ satisfying
\begin{align}\label{eq:def operator delbat J,H}
\del_s u + J(u)(\del_t u - X_H(u))=0
\end{align}
and which admit a lift $\tilde{u}\co \bb R\ra \widetilde{\m LM}$ with limits $c$ and $c'$. Indeed, $\bb R$ acts on $\m M(c,c';H,J)$ by translation and the regularity condition satisfied by the pair $(H,J)$ ensures that, for $c$ and $c'$ with index difference $1$, $\m M(c,c';H,J)/\bb R$ is a compact $0$--dimensional manifold.

Floer homology is the homology of this complex, and does not depend, up to natural isomorphism, on the choice of the regular pair $(H,J)$: $HF_*(M,\omega)=H_*(CF(H),\del_{(H,J)})$. This natural isomorphism, \defin{the comparison morphism}, appears explicitly in the construction of Seidel's representation.\\

Indeed, Seidel's morphism can be seen as a representation of $\tilde{\pi}_1(\ham{M,\omega})$, namely, for any $[g,\tilde{g}]\in\tilde{\pi}_1(\ham{M,\omega})$, one can associate an automorphism
\begin{align*}
\seidel[g,\tilde{g}] \co HF_*(M,\omega)\lra HF_{*-2I(g,\tilde{g})}(M,\omega)
\end{align*}
defined as the composition $\seidel[g,\tilde{g}]=H_*(\comp({\bf H},{\bf J})\circ\nat(g,\tilde{g}))$. 

The morphism $\comp({\bf H},{\bf J})$ is the comparison morphism of Floer homology
\begin{align*}
\comp({\bf H},{\bf J}) \co CF_*(M,\omega;H_0,J_0)\lra CF_*(M,\omega;H_1,J_1).
\end{align*}
It is defined by using $({\bf H},{\bf J})$, any regular homotopy between $(H_0,J_0)$ and $(H_1,J_1)$ and it induces an isomorphism in homology (see for example \cite{Salamon99} for the definition and the isomorphism and naturality properties of this morphism).

The second morphism involved in the composition defining $\seidel[g,\tilde{g}]$ is \defin{the naturality morphism} which is an identification of chain complexes, coming from the action of $(g,\tilde{g})$ on $\widetilde{\m LM}$,
\begin{align*}
\nat(g,\tilde{g})\co CF_*(M,\omega;H^g,J^g)\lra CF_{*-2I(g,\tilde{g})}(M,\omega;H,J)
\end{align*}
where the pair $(H^g,J^g)$ is the pushforward of $(H,J)$ by $g$, defined as
\begin{align}\label{eq:pushforward of pairs}
\begin{split}
H^g(t,y) &= H(t,g_t(y))-K_g(t,g_t(y)),\\
\mathrm{and}\; J^g_t &= dg_t^{-1}J_t dg_t
\end{split}
\end{align}
($K_g$ being a Hamiltonian generating the loop $g$). By straightforward computations, one can show that $(H^g,J^g)$ is a regular pair if and only if $(H,J)$ is regular.

Now, to define $\seidel[g,\tilde{g}]$ as the composition of these two morphisms, we choose the pair $({\bf H},{\bf J})$, required to define $\comp({\bf H},{\bf J})$, to be any regular homotopy between $(H_0,J_0)=(H,J)$ and $(H_1,J_1)=(H^g,J^g)$.

The definition of the shift of indices $I(g,\tilde{g})$ is standard (it corresponds to the degree of a loop in $\mathrm{Sp}(2n,\bb R)$ coming from a trivialization of $TM$ over the cappings $v$'s of the orbits $x$'s -- see the definition of $\widetilde{\m LM}$). It is compatible with the Cartesian product, in the following sense:
\begin{align*}
I(g,g';\tilde{g},\tilde{g}')= I(g,\tilde{g}) + I(g',\tilde{g}').
\end{align*}
In view of this formula, the shift of indices will be implied in what follows.\\

The correspondence between the two descriptions of Seidel's representation, namely, between $q[g,\tilde{g}]$ and $\seidel[g,\tilde{g}]$, goes via the Piunikhin--Salamon--Schwarz (PSS) morphism as well as the pair-of-pants product. These tools appeared in \cite{PiunikhinSalamonSchwarz96}. We recall that (under the monotonicity assumption) the PSS morphism is a canonical isomorphism
\begin{align*}
\PSS\co QH_*(M,\omega)\lra HF_*(M,\omega)
\end{align*}
between the quantum homology and the Floer homology of $(M,\omega)$, as modules over the Novikov ring. The pair-of-pants product is a product on Floer homology
\begin{align*}
\pop \co HF_*(M,\omega)\otimes HF_*(M,\omega)\lra HF_*(M,\omega)
\end{align*}
defined on chain complexes by counting suitable moduli spaces of pair-of-pants. Given these tools, Seidel proved that
\begin{align}\label{eq:link between q and S}
\seidel[g,\tilde{g}](b)=\PSS(q[g,\tilde{g}])\pop b
\end{align}
for all $b\in HF_*(M,\omega)$. This is the interpretation we use to prove Theorem \ref{theo:main theo}.


\section{Proof of the theorem}\label{sec:proof of theo}

Let $(M,\omega)$ and $(M',\omega')$ be closed monotone symplectic manifolds (with the same monotonicity constant). Let $(H,J)$ and $(H',J')$ be respectively defined on $(M,\omega)$ and $(M',\omega')$. We define on $M\times M'$ the Hamiltonian $H\oplus H'$ and the almost complex structure $J\oplus J'$ by the formulae
\begin{align}\label{eq:defin split pairs}
\begin{split}
(H\oplus H')_t(x,x') &=H_t(x) + H_t(x')\\
(J\oplus J')_t(\xi,\xi') &=(J_t(\xi),J'_t(\xi'))
\end{split}
\end{align}
for all $(x,x')\in M\times M'$ and all $(\xi,\xi')\in T_xM\otimes T_{x'}M' \simeq T_{(x,x')}(M\times M')$.

\begin{rema}\label{rema:pushforward product is product of pushforwards}
Notice that the pushforward, as defined by (\ref{eq:pushforward of pairs}), of $(H\oplus H',J\oplus J')$ by any element of the form $(g,g')\in\ham{M\times M',\omega\oplus\omega'}$ satisfies
    \begin{align*}
    \big((H\oplus H')^{(g,g')}, (J\oplus J')^{(g,g')}\big)=(H^g\oplus H'^{g'}, J^g\oplus J'^{g'}).
    \end{align*}
\end{rema}

Even for regular pairs $(H,J)$ and $(H',J')$, the pair $(H\oplus H',J\oplus J')$ is not \textit{a priori} regular. As we shall see, the problem comes from the fact that the moduli spaces of simple spheres of the product is, in general, bigger than the product of the moduli spaces of simple spheres of each component. Thus, the complex structure $J\oplus J'$ is not automatically regular. In \S \ref{section:regularity of split pairs}, we show that to go through the construction a weaker regularity condition is enough. This will give sense to the following claim.

\begin{claim}\label{claim:H,J regular enough}
If $(H,J)$ and $(H',J')$ are regular pairs, the pair $(H\oplus H',J\oplus J')$ is regular enough. Moreover, if $(\mathbf{H},\mathbf{J})$ and $(\mathbf{H}',\mathbf{J}')$ are regular homotopies, then the homotopy $(\mathbf{H}\oplus\mathbf{H}',\mathbf{J}\oplus\mathbf{J}')$ is regular enough.
\end{claim}
We postpone the proof until the next section. As mentioned in Remark \ref{rema:on the monotonicity restriction}, the proof of the part concerning almost complex structures is the only place where we use (a property implied by) the monotonicity assumption.\\

Now, notice that
\begin{align*}
\mathrm{Crit}_k(\m A_{H\oplus H'}) = \bigcup_{l+l'=k} \mathrm{Crit}_l(\m A_H) \times \mathrm{Crit}_{l'}(\m A_{H'})
\end{align*}
and that the action agrees with this decomposition, that is, for all $[v,x]\in \mathrm{Crit}_l(\m A_H)$ and all $[v',x']\in \mathrm{Crit}_{l'}(\m A_{H'})$:
\begin{align*}
\m A_{H\oplus H'}([v,v';x,x'])=\m A_{H}([v,x]) + \m A_{H'}([v',x']).
\end{align*}
The finiteness condition (\ref{eq:finiteness condition FH}) is such that
\begin{align*}
\bigoplus_{l+l'=k} CF_l(H)\otimes CF_{l'}(H') \simeq CF_k(H\oplus H').
\end{align*}
This isomorphism induces a morphism in homology
\begin{align*}
\kappa\co HF_*(M,\omega)\otimes HF_*(M',\omega') \lra HF_*(M\times M',\omega\oplus\omega').
\end{align*}

\begin{claim}\label{claim:diag commutes 1}
  The following diagram commutes
\begin{align*}
  \xymatrix{\relax
        HF_*(M,\omega)\otimes HF_*(M',\omega') \ar[r]^{\kappa}\ar[d]_{\seidel[g,\tilde{g}]\otimes \seidel[g',\tilde{g'}]} & HF_*(M\times M',\omega\oplus\omega')\ar[d]^{\seidel[g,g';\tilde{g},\tilde{g}']}\\
        HF_*(M,\omega)\otimes HF_*(M',\omega') \ar[r]^{\kappa} & HF_*(M\times M',\omega\oplus\omega')
      }
\end{align*}
\end{claim}

\begin{proof}[Proof of Claim \ref{claim:diag commutes 1}]
  Decomposing the automorphisms of Floer homology given by $[g,\tilde{g}]$, $[g',\tilde{g}']$ and $[g,g';\tilde{g},\tilde{g}']$ in terms of naturality and comparison morphisms, we want to prove that the following diagram commutes in homology
\begin{align*}
  \xymatrix{\relax
        CF_*(H^g,J^g)\otimes CF_*(H'^{g'},J'^{g'}) \ar[r]\ar[d]_{\nat[g,\tilde{g}]\otimes \nat[g',\tilde{g'}]} & CF_*(H^g\oplus H'^{g'}, J^g\oplus J'^{g'})\ar[d]^{\nat[g,g';\tilde{g},\tilde{g}']}\\
        CF_*(H,J)\otimes CF_*(H',J') \ar[r]\ar[d]_{\comp(\mathbf{H},\mathbf{J})\otimes \comp(\mathbf{H}',\mathbf{J}')} & CF_*(H\oplus H', J\oplus J') \ar[d]^{\comp(\mathbf{H}\oplus\mathbf{H}',\mathbf{J}\oplus\mathbf{J}')}\\
        CF_*(H^g,J^g)\otimes CF_*(H'^{g'},J'^{g'}) \ar[r] & CF_*(H^g\oplus H'^{g'}, J^g\oplus J'^{g'})
      }
\end{align*}
where $(H,J)$, $(H',J')$, etc are defined as above. The diagram even commutes at the chain level, with the choices we made (justified by Claim \ref{claim:H,J regular enough}) and by Remark \ref{rema:pushforward product is product of pushforwards}: The horizontal maps identify products of moduli spaces with moduli spaces of the product (for any type of moduli spaces involved by these morphisms).
\end{proof}

By (\ref{eq:link between q and S}) which relates the two descriptions of Seidel's morphism, Claim \ref{claim:diag commutes 1} immediately amounts to the fact that for $b\in HF_*(M,\omega)$ and $b'\in HF_*(M',\omega')$,
\begin{align}\label{eq:PSS o i = i o PSS}
\PSS(q[g,g';\tilde{g},\tilde{g}'])\pop \kappa(b\otimes b') = \kappa\big( (\PSS(q[g,\tilde{g}])\pop b)\otimes (\PSS(q[g',\tilde{g'}])\pop b')\big).
\end{align}

\begin{claim}\label{claim:diag commutes 2}
The following diagram commutes
\begin{align*}
\xymatrix{\relax
  QH_*(M,\omega)\otimes QH_*(M',\omega') \ar[r]^{\kappa_Q}\ar[d]_{\PSS\otimes \PSS} & QH_*(M\times M',\omega\oplus\omega') \ar[d]^{\PSS}\\
  HF_*(M,\omega)\otimes HF_*(M',\omega') \ar[r]^{\kappa} & HF_*(M\times M',\omega\oplus\omega')
}
\end{align*}
\end{claim}

\begin{proof}[Proof of Claim \ref{claim:diag commutes 2}]
If the parameters (Hamiltonian functions, almost complex structures, Morse functions, metrics, etc) used to define the involved PSS morphisms are chosen as above, namely split, the products of moduli spaces are again identified with the moduli spaces of the product and the commutativity even holds at the chain level.
\end{proof}

For the element $q[g,\tilde{g}]\otimes q[g',\tilde{g}']$, this commutativity amounts to
\begin{align}\label{eq:PSS o i = i o PSS -- 2}
  \kappa\big(\PSS(q[g,\tilde{g}])\otimes \PSS(q[g',\tilde{g}'])\big) = \PSS\big(\kappa_Q(q[g,\tilde{g}]\otimes q[g',\tilde{g'}])\big).
\end{align}
Since $[M\times M'] = \kappa_Q([M]\otimes [M'])$, Claim \ref{claim:diag commutes 2} also leads to
\begin{align*}
\PSS([M\times M']) &= \PSS\big(\kappa_Q([M]\otimes [M'])\big) = \kappa\big(\PSS([M])\otimes \PSS([M'])\big).
\end{align*}

Thus, letting $b=\PSS([M])$ and $b'=\PSS([M'])$ in (\ref{eq:PSS o i = i o PSS}), one gets that
\begin{align}\label{eq:PSS o i = i o PSS -- 3}
\PSS(q[g,g';\tilde{g},\tilde{g}']) = \kappa\big(\PSS(q[g,\tilde{g}])\otimes \PSS(q[g',\tilde{g'}])\big)
\end{align}
since the image via the PSS morphism of the fundamental class (the identity element of the group of invertible elements of quantum homology) acts trivially for the pair-of-pants product. Finally, (\ref{eq:PSS o i = i o PSS -- 2}) and (\ref{eq:PSS o i = i o PSS -- 3}) amount to
\begin{align*}
  \PSS\big(q[g,g';\tilde{g},\tilde{g}']\big) = \PSS\big(\kappa_Q(q[g,\tilde{g}]\otimes q[g',\tilde{g'}])\big).
\end{align*}
This completes the proof of the theorem, since the PSS morphism is an isomorphism.


\section{Regularity of split pairs}\label{section:regularity of split pairs}

In this section, we give precise definitions of regularity (for almost complex structures and pairs $(H,J)$), we define ``regular enough pairs'' and prove Claim \ref{claim:H,J regular enough}. We consider the case of $S^1$--families of $\omega$--compatible almost complex structures. This is sufficient to prove that Floer homology is well-defined. The case of $2$--parameter families of almost complex structures (needed for instance for homotopies) works along the same lines.\\

Let $\m M^s(J)$ denote the set of pairs $(t,w)\in S^1\times C^\infty(S^2,M)$, where $w$ is a $J_t$--holomorphic \textit{simple} sphere in $M$. This set is the union over $k$ of the subsets $\m M^s_k(J)$ of pairs with spheres of first Chern number $k$. With $S^1$--families of almost complex structures, the linearization of the equation $\bar{\del}_J=0$ at $(t,w)$ is given by
\begin{align}\label{eq:def linear operator for J}
\begin{split}
\hat{D}_J(t,w)\co T_t S^1\times C^\infty(w^*TM)\lra \Omega^{0,1}(w^*(TM,J_t))\\
\hat{D}_J(t,w)(\theta,W)= D_{J_t}(w)W + \frac{1}{2} DJ(t)\theta\circ dw\circ i
\end{split}
\end{align}
where $i$ is the complex structure of $S^2\simeq \bb C\bb P^1$ and $DJ(t)$ denotes the derivative of the $S^1$--family of almost complex structures at $t$.

\begin{defi}\label{def:regularity condition on acs}
The $S^1$--family of almost complex structures $J$ is \defin{regular} if the linearized operator defined by (\ref{eq:def linear operator for J}) is onto for all $(t,w)\in \m M^s(J)$.
\end{defi}

Now denote by $V_k(J)\subset S^1\times M$ the set of pairs $(t,x)$ for which there exists a non-constant, $J_t$--pseudo-holomorphic sphere $w$, with first Chern number $c_1(w)\leq k$ and such that $x\in \im(w)$.

\begin{defi}\label{def:regularity condition on pairs H,J}
A pair $(H,J)$ consisting of a family of almost complex structures $J$ and a time-dependent Hamiltonian $H\co S^1\times M\ra \bb R$ is \defin{regular} if
\begin{itemize}
\item[\textbf{i.}] $J$ is a regular $S^1$--family of almost complex structures,
\item[\textbf{ii.}] the critical points of $\m A_H$ are non-degenerate and for any orbit $x$ of the Hamiltonian vector field and all $t\in S^1$, $(t,x(t))\notin V_1(J)$,
\item[\textbf{iii.}] the linearization of the operator (\ref{eq:def operator delbat J,H}) is onto for all $u \in \m M(c,c';H,J)$, and
\item[\textbf{iv.}] if $\mathrm{ind}(u)\leq 2$ then for all $t$ and $s$, $(t,u(s,t))\notin V_0(J)$.
\end{itemize}
\end{defi}

It is well-known that for monotone symplectic manifolds the sets of regular almost complex structures and of regular pairs are dense.

\subsection{Regular enough almost complex structures and pairs}

When $J$ is regular, the following fundamental claims hold.

\begin{claim2}\label{claim:MkJ smooth mfds}
For all integers $k$, the set $\m M_k^s(J)$ is (either empty or) a smooth manifold of dimension $2n+2k+1$.
\end{claim2}

\begin{claim2}\label{claim:no hol S2 of negative chern}
$J$ is semi-positive, that is, for all $k<0$, $\m M_k^s(J)$ is empty.
\end{claim2}

We can now define a weaker regularity condition for almost complex structures.

\begin{defi}
An almost complex structure is \defin{regular enough} if Claim \ref{claim:MkJ smooth mfds} holds for $0\leq k\leq 2$ and Claim \ref{claim:no hol S2 of negative chern} holds. A pair is \defin{regular enough} if the almost complex structure is regular enough and the pair satisfies conditions \textbf{ii}-\textbf{iv} of Definition \ref{def:regularity condition on pairs H,J}.
\end{defi}

Roughly speaking, Claim \ref{claim:MkJ smooth mfds}, for $k=0$ and $1$ implies that $V_0(J)$ has codimension $4$ and $V_1(J)$ has codimension $2$ as subsets of $S^1\times M$. Thus, for a regular enough almost complex structure, the choice of a Hamiltonian $H$ such that the pair $(H,J)$ satisfies conditions \textbf{ii} and \textbf{iv} of Definition \ref{def:regularity condition on pairs H,J} is generic.

Now, for such a pair, Floer homology is well-defined since bubbling is avoided. Indeed, condition \textbf{ii} forbids configurations of the type index--$0$ tube (that is, a ``constant'' tube which coincides with an orbit) with attached spheres whose first Chern numbers sum to $1$. Condition \textbf{iv} forbids the appearance of configurations of the type index--$2$ tube with attached spheres whose first Chern numbers sum to $0$.

\subsection{Split pairs are regular enough (proof of Claim \ref{claim:H,J regular enough})}

Let $(M,\omega)$ and $(M',\omega')$ be closed monotone symplectic manifolds (with the same monotonicity constant). A \defin{split pair} on $M\times M'$ is a pair of the form $(H\oplus H',J\oplus J')$, as defined by (\ref{eq:defin split pairs}), where $(H,J)$ and $(H',J')$ are respectively defined on $(M,\omega)$ and $(M',\omega')$.\\

If $J$ and $J'$ are compatible with the respective symplectic forms, $J\oplus J'$ is obviously $(\omega\oplus\omega')$--compatible.

\begin{lemm}\label{lemm:regularity of J+J'}
If $J$ and $J'$ are regular, then $J\oplus J'$ is regular enough.
\end{lemm}

The main observation is the decomposition of $\m M^s_k(J\oplus J')$ in terms of pseudo-holomorphic spheres in $M$ and $M'$. 
More precisely, a simple sphere $w=(v,v')$ of $M\times M'$ is a pair of
\begin{itemize}
\item[$-$] two simple spheres,
\item[$-$] a simple sphere and a constant sphere,
\item[$-$] a simple sphere and a multiply covered sphere,
\item[$-$] two multiply covered spheres (with relatively prime degrees).
\end{itemize}
Notice that a pair consisting of a constant sphere and a multiply covered one is a multiply covered sphere (of the product).

\begin{rema}
Since the linearized operator defined by (\ref{eq:def linear operator for J}) respects the splitting $w^*T(M\times M')=v^*TM\times v'^*TM'$, it is onto at $(t,w)$ if and only if both its projections are onto (at $(t,v)$ and $(t,v')$). Thus a split almost complex structure is not \textit{a priori} regular, since we have no information about the operator related to $J$ (for instance) at $(t,v)$ for multiply covered $v$.
\end{rema}

\begin{proof}[Proof of Lemma \ref{lemm:regularity of J+J'}]
We denote by $\m M^{m,i}_l(J)$ the set of pairs $(t,v)$ where $v$ is an $i$--fold covered, $J_t$--pseudo-holomorphic sphere in $M$ with first Chern number $l$. Notice that (by convention) $\m M^{m,1}_l(J)=\m M^s_l(J)$, but $\m M^{m}_l(J)$ denotes the union of these sets over $i>1$ (that is, $\m M^{m}_l(J)$ is the set of pairs whose sphere is \textit{strictly} multiply covered). We get, for any integer $k$,
\begin{align}
\m M^s_k(J\oplus J') \simeq
&      && \bigcup_{\;\;\;l+l'=k\;\;\;}  && \m M^s_l(J) \times_{S^1} \m M^s_{l'}(J') & \label{simple} \\
& \cup &&                               && \m M^s_k(J)  \times M' \;\cup\; M \times \m M^s_{k}(J') & \label{simple and constant} \\
&\cup&&\bigcup_{\;\;\;l+l'=k\;\;\;}&& \m M^s_l(J)\times_{S^1} \m M^m_{l'}(J') \;\cup\; \m M^m_l(J) \times_{S^1} \m M^s_{l'}(J') & \label{simple and multiply} \\
& \cup && \bigcup_{\substack{l+l'=k,\\ i,j \geq 2\\ \mathrm{gcd}(i,j)=1}} && \m M^{m,i}_l(J) \times_{S^1} \m M^{m,j}_{l'}(J') & \label{multiply}
\end{align}
(we recall that the fibered products $\m M(J) \times_{S^1} \m M(J')$ appearing above consist of the elements $\big((t,v),(t',v')\big)$ of $\m M(J) \times \m M(J')$ for which $t=t'$). The union of (\ref{simple}), (\ref{simple and multiply}) and (\ref{multiply}) can be described as
\begin{align*}
\bigcup_{\substack{l+l'=k,\\ i,j \geq 1|\,\mathrm{gcd}(i,j)=1}} \m M^{m,i}_l(J) \times_{S^1} \m M^{m,j}_{l'}(J').
\end{align*}
However, (\ref{simple}) and (\ref{simple and constant}) are really different from (\ref{simple and multiply}) and (\ref{multiply}) and should be studied separately.\\

Since $J$ and $J'$ are regular, there is no pseudo-holomorphic sphere with negative first Chern number (Claim \ref{claim:no hol S2 of negative chern}). Thus, from the decomposition above, we can already conclude that Claim \ref{claim:no hol S2 of negative chern} holds for $J\oplus J'$. Moreover, for a non-empty set appearing in the decomposition, we have $l$ and $l'$ non-negative and furthermore if $k=0$, then $l=l'=0$. Let now look at small values of $k$.\\

Since a non-constant pseudo-holomorphic sphere has positive symplectic area, such a sphere in a monotone symplectic manifold cannot have a vanishing first Chern number. Hence $\m M^s_0(J\oplus J') = \emptyset$. Moreover, the first Chern number of a non-constant pseudo-holomorphic \textit{(strictly) multiply covered} sphere is at least $2$. Thus,
\begin{align*}
\m M^s_1(J\oplus J') = \m M^s_1(J)  \times M' \;\cup\; M \times \m M^s_{1}(J')
\end{align*}
which is a smooth manifold of the expected dimension $2(n+n')+3$ (see Claim \ref{claim:MkJ smooth mfds}). Finally for $k=2$, in the decomposition above we have either $l=1$ and $l'=1$ (and there is no multiply covered pseudo-holomorphic sphere), or $l=0$ or $l'=0$ (and then at least one sphere has to be constant). Since a pair consisting of a constant sphere and a multiply covered one is not simple, we can conclude that
\begin{align*}
\m M^s_2(J\oplus J') = \m M^s_1(J) \times_{S^1} \m M^s_{1}(J') \;\cup\; \m M^s_2(J)  \times M' \;\cup\; M \times \m M^s_2(J')
\end{align*}
which is the union of three smooth manifolds, of the expected dimension, $2(n+n')+5$ (see Claim \ref{claim:MkJ smooth mfds}). This proves that Claim \ref{claim:MkJ smooth mfds} holds when $k=0$, $1$ and $2$.
\end{proof}

\begin{rema}\label{rema:SSP non monotone}
Lemma \ref{lemm:regularity of J+J'} is the only place where we \textit{a priori} have to restrict our study to monotone symplectic manifolds. The main problem appearing in the general (strongly semi-positive) case, comes from the existence of non-constant pseudo-holomorphic spheres with vanishing first Chern number. When such spheres exist, the moduli spaces $M^s_k(J\oplus J')$ (for $k=0$, $1$ and $2$) are more complicated.

However, the additional subsets are products of smooth manifolds (moduli spaces of constant or simple pseudo-holomorphic spheres for regular almost complex structures) and of moduli spaces of strictly multiply covered spheres: $\m M^{m,i}_k(J)$ and $\m M^{m,j}_k(J')$ with $i$, $j>1$ and $k=0$ or $2$. Such moduli spaces are formed of finitely many copies of sets, in bijection with smooth manifolds of the expected dimension (or of codimension at least $2$ in the whole union).

Indeed, for a $i$--fold covered sphere to represent the homotopy class $[A]\in \pi_2(M)$, there has to be a primitive homotopy class $[B]$, with $[A]=i[B]$. Thus, for each homotopy class $[A]$, there are only finitely many integers $i$ for which $[A]$ admits a representative which is $i$--fold covered. (Since we consider non-trivial pseudo-holomorphic spheres, they represent non-zero homotopy classes).

Now, for each of these integers, the set $\m M^{m,i}_k(J)$ is in bijection with $i$ disjoint copies of $\m M^{s}_{\frac{k}{i}}(J)$ (by considering the underlying simple curves, the canonical degree--$i$ map of $\bb C\bb P^1$: $[z_1:z_2]\mapsto [z_1^i:z_2^i]$ and the action of each $i$--th root of unity $\xi$: $[z_1:z_2]\mapsto [\xi z_1:z_2]$). By regularity of $J$, $\m M^{s}_\frac{k}{i}(J)$ is of dimension: $\mathrm{dim}(M) + 2\frac{k}{i}$, that is, the expected dimension when $k=0$ and the expected dimension minus (at least) two when $k>0$.
\end{rema}

\begin{lemm}\label{lemm:regularity of split pairs}
If the pairs $(H,J)$ and $(H',J')$ are regular, then $(H\oplus H', J\oplus J')$ is regular enough.
\end{lemm}

\begin{proof}
Lemma \ref{lemm:regularity of J+J'} proves that for such pairs $J\oplus J'$ is regular enough. We now prove conditions \textbf{ii}-\textbf{iv} of Definition \ref{def:regularity condition on pairs H,J}.

\hop [\textbf{ii}.] -- By definition of the split Hamiltonian,
\begin{align*}
\mathrm{Crit}(\m A_{H\oplus H'})=\mathrm{Crit}(\m A_{H})\times \mathrm{Crit}(\m A_{H'}).
\end{align*}
Thus, the critical points of $\m A_{H\oplus H'}$ are non-degenerate.

For such a critical point $[v,v';x,x']$, and some $t\in S^1$, $(t;x(t),x'(t))\in V_1(J\oplus J')$ amounts to the existence of a non-constant $(J\oplus J')_t$--pseudo-holomorphic sphere $w\co S^2\ra M\times M'$, with first Chern number less or equal to $1$, and passing through $(x(t),x'(t))$. We recall that there is no pseudo-holomorphic sphere (in $M$ and $M'$, for $J$ and $J'$) with negative first Chern number. Thus, $w$ is the product of pseudo-holomorphic spheres in $M$ and $M'$, one of them (let say the component in $M$) being non-constant, with first Chern number $0$ or $1$.
Thus, $(t,x(t))$ lies in $V_1(J)$ which contradicts the fact that $(H,J)$ satisfies \textbf{ii}. Hence, such a $w$ does not exist and $(H\oplus H',J\oplus J')$ satisfies \textbf{ii}.

\hop [\textbf{iii}.] -- The moduli spaces defining the differential of the Floer complex of the product split, that is,
\begin{align*}
\m M\big((c,c'),(d,d');H\oplus H',J\oplus J'\big)=\m M(c,d;H,J)\times \m M(c',d';H',J')
\end{align*}
and thus so does the operator $\bar\del_{J\oplus J',H\oplus H'}$. Thus the linearized operator is onto if and only if both operators (in each component of the product) are onto; since $(H,J)$ and $(H',J')$ satisfy \textbf{iii}, so does $(H\oplus H',J\oplus J')$.

\hop [\textbf{iv}.] -- In view of the decomposition of the moduli spaces above, if the index of $(u,u')\in \m M(c,d;H,J)\times \m M(c',d';H',J')$ is less or equal than $2$, then $\mathrm{ind}(u)\leq 2$ and $\mathrm{ind}(u')\leq 2$. Then we can conclude, as in the proof of \textbf{ii}, that the existence of some $t$ such that $(t;u(s,t),u'(s,t))\in V_0(J\oplus J')$ implies that $(t,u(s,t))\in V_0(J)$ (and/)or $(t,u'(s,t))\in V_0(J')$ . Thus, since $(H,J)$ and $(H',J')$ satisfy \textbf{iv}, so does $(H\oplus H',J\oplus J')$. (Notice that, in the monotone case, the set of non-constant pseudo-holomorphic spheres of first Chern number $0$ in the product is empty and thus that condition \textbf{iv} is trivially satisfied. The previous justification does not use the assumption of monotonicity.)
\end{proof}

Lemma \ref{lemm:regularity of split pairs} proves Claim \ref{claim:H,J regular enough} for pairs consisting of time-dependent Hamiltonians and almost complex structures. The same arguments can be carried out for homotopies.

\end{document}